\documentclass{article}

\newcommand{\C}{I\!\! C}

\begin{document}

Tsemo Aristide

  University of Toronto

   St. George Street

Toronto, Ontario Canada

M5S 3G3

tsemo58@yahoo.ca

\bigskip

\centerline{\bf The structure of CR manifolds.}

\bigskip

\centerline{Abstract.}

{\it In this paper we study the topology of CR pseudoconvex
manifolds whose Reeb flow preserve the Levi metric.}

\bigskip
\bigskip

{\bf Definition 1.}

A CR-manifold of dimension $2n+1$ is a manifold $N$ of dimension
$2n+1$ endowed with the followings properties:

Let $TN$ be the tangent bundle of $N$, there exists a subbundle
$V$ of $TN\otimes{\C}$, of complex dimension $n$.  such that
$V\cap \bar V=0$, where $\bar V$ is the complex conjugate of $V$,
and $[V,V]\subset V$.

The bundle $V\oplus \bar V=H$ is the kernel of a real $1$-form
$\theta\in T^*N$, where $T^*N$ is the cotangent bundle of $N$.

Let $(U_p)_{p\in I}$ be an open contractible covering of $N$ by
open subsets. We denote by $\theta_p$ the restriction of $\theta$
to $U_p$, and $(u^p_1,...,u^p_n)$ a basis of the restriction of
the dual of $V$ to $U$, we have:

$$
d(\theta_p)=\sum ih_{cd}u_c\wedge \bar{u_d}
$$

The coefficients $h_{cd}$, define a pseudo-Hermitian metric on
$V$. We can extend this metric to $H$ by supposing that $V$ is
orthogonal to $\bar V$, and the complex conjugate is an isometry.
This metric is called the Levi metric.

\medskip

{\bf Definition 2.}

The CR-manifold is pseudo convex if the Levi metric is positive
definite. In this case there exists a vector field $X$ of $TN$
such that $L_X\theta=0$, where $L_X$ is the Lie derivative. The
vector field $X$ is called the Reeb vector field.

The purpose of this paper is to show the following result:

\bigskip

{\bf Theorem 1.}

{\it Let $N$ be a compact pseudo-convex CR manifold, $X$ a flow
transverse to $H$ which preserves the Levi metric then the closure
of the orbits of  $X$ are torus, if these closure have the same
dimension, then $N$ is the total space of a bundle whose typical
fiber is the closure of an orbit of $X$ flow over a Satake
manifold.}

\bigskip

{\bf The complex transverse geometry of the Reeb flow.}

\bigskip

Let $N$ be a CR pseudo-convex manifold, and $X$ a flow transverse
to $H$ which preserves the Levi metric. The complex transverse
bundle $T_XN$ of $X$ is the quotient of $TN\otimes{\C}$  by
$X\otimes{\C}$. Let $U$ be a contractible open subset of $N$, and
$T$ a submanifold of $U$ transverse to $X$. We denote by
${L_XN}_U$, the restriction of the bundle of complex frames $L_XN$
of $T_XN$ to $U$, we have a projection $p_U:L_XN\rightarrow
LT\otimes {\C}$, where $LT$ is the bundle of transverse complex
frames  of $T$. The kernel of the differential $dp_U$ of $d_P$
define on ${L_XN}_U$  a distribution tangent to a flow. This
distribution is independent of the choice of the local transversal
$T$. We have thus defined on $L_XN$ a flow $\hat X$, which is
called the lift of $X$.

The Levi metric $<,>$, defines an Hermitian reduction $H_X$ of
$L_XN$ (which is invariant by the orbit of $\hat X$) of $L_XN$,
since the manifold $N$ is compact, $H_X$ is also a compact
submanifold. To show the theorem 1, we shall prove that the
restriction $X'$ of $\hat X$ to $H_X$ is a riemaniann flow, (in
fact, we are going to construct  a transverse parallelism to $X'$)
and apply a well-known result of Yves Carriere on the structure of
Riemannian flows.

\bigskip

{\bf The transverse parallelism of $X'$.}

\bigskip

The manifold $H_X$ is a principal bundle over $N$ whose typical
fiber is $U(2n)$, the group of Hermitian matrices. Each element
$c$ of the Lie algebra  $u(2n)$ of $U(2n)$ defines on $H_X$ a
vector field $c_N$ defined by
$c_N(x)={d\over{dt}}_{t=0}exp(tc_N)(x)$. The vector $c_N$ are
called the fundamental vector fields of $H_X$

\medskip

Let $\alpha$ be the fundamental form of $H_X$. Recall that for
every element $u$ of the tangent space ${TH_X}_x$ of $x\in TH_x$
($x$ is a linear map ${\C}^{2n}\rightarrow T_{p(x)}N\otimes
{\C}/{\C}\otimes X$, where $p:H_X\rightarrow N$ is the bundle
projection map, $\alpha(u)=x^{-1}(dp_x(u))$. For each element $y$
of ${\C}^{2n}$, we can define the vector field $\hat y$ of $H_X$
by setting $\alpha(\hat y)=y$.

The vector $\hat y$, $y\in {\C}^{2n}$ and the fundamental vector
fields define the transverse parallelism to $X'$.

\bigskip

{\bf Proof of theorem 1.}

\medskip

Let $X$ be an isometry transverse to the Levi metric, then the
lift $X'$ of $X$ to $H_X$ is a riemannian flow, we can apply the
result of Carriere. The closure of the orbits of $X$ are the
projections of the closure of the orbits of $X'$ by the bundle map
$H_X\rightarrow N$.

\bigskip

We can obtain this most general result:

\bigskip

{\bf Theorem 2.}

{\it Let $N$ be a pseudo convex compact CR manifold, endowed with
a transverse flow $X$  which is a conformal flow in respect to the
Levi metric, then the closure of the orbits of $X$ are torus.}

\medskip

{\bf Proof.}

  Consider the group generated by $U(n)$ and the complex
homothetic maps, and we denote $L(n)$, its quotient by an
homothetic map $h_{\lambda_0}$ such that the norm of $\lambda_0$
is strictly superior to $1$. The transverse bundle of $X$ can be
reduced to $L(n)$. The lifts $X'$ of $X$ to $L(n)$ is a
transverselly conformal analytic flow of codimension greater than
$3$. We can apply the theorem of Tarquini which asserts that in
this situation that the flow of $X'$ is riemannian or Moebius.
Then we apply the structure theorem for riemannian flows and
Moebius flows.

\bigskip

The previous result suggests the study of transversely Hermitian
foliations:

\medskip

{\bf Definition 3.}

Let $N$ be a manifold, and ${\cal H}$ a foliation defined on $N$,
the foliation ${\cal H}$ is a transversely Hermitian foliation, if
there exists a symmetric Hermitian basic $2$-tensor $<,>$ defined
on $TN\otimes {\C}$ such that:

For each $u\in N$, $T_u{\cal H}\otimes {\C}$ the tensor product of
the subspace $T_x{\cal H}$ of $T_uN$ tangent to ${\cal H}$ and
${\C}$ is the kernel of $<,>$.

The projection of $<,>$ to the transverse complex bundle of ${\cal
H}$ is a positive definite Hermitian metric. We have the following
nice result:

\bigskip

{\bf Theorem 3.}

{\it Let $N$ be a compact manifold endowed with a transversely
Hermitian foliation ${\cal H}$, then the closure of the leaves of
${\cal H}$ are submanifolds.}

\bigskip

\centerline{\bf References.}

1. Yves Carriere, Flots riemanniens Asterisque 116

2. Tarquini, Feuilletages conformes, Annales Institut Fourier 2004

\end{document}